\documentclass{amsart}

\usepackage[utf8x]{inputenc}
\usepackage{bbm}
\usepackage{yfonts}
\usepackage{amsfonts}
\usepackage{amsmath}
\usepackage{amssymb}
\usepackage{amsthm}
\usepackage{amscd}
\usepackage{tensor}
\usepackage{graphicx}
\usepackage[all,cmtip]{xy}
\usepackage{stmaryrd}
\usepackage{tikz-cd}
\usepackage{tikz}
\usetikzlibrary{matrix,arrows,decorations.pathmorphing}
\usepackage{calc}
\usepackage[cal=boondox]{mathalfa}
\usepackage[toc,page]{appendix}
\usepackage{multicol}

\usepackage[hyperfootnotes=false]{hyperref}
\hypersetup{colorlinks=true,pageanchor=false,linkcolor=blue,citecolor=red}
\usepackage[all]{hypcap}

\begin{document}

\pdfsuppresswarningpagegroup=1

\newtheorem{theorem}{Theorem}[section]
\newtheorem{corollary}[theorem]{Corollary}
\newtheorem{proposition}[theorem]{Proposition}
\newtheorem{lemma}[theorem]{Lemma}

\theoremstyle{definition}
\newtheorem{definition}[theorem]{Definition}
\newtheorem{example}[theorem]{Example}
\newtheorem{condition}[theorem]{Condition}
\newtheorem{construction}[theorem]{Construction}
\newtheorem{exercise}[theorem]{Exercise}

\theoremstyle{remark}
\newtheorem{remark}[theorem]{Remark}
\newtheorem{fact}[theorem]{Fact}
\newtheorem{claim}[theorem]{Claim}
\newtheorem{reminder}[theorem]{Reminder}
\newtheorem{conjecture}{Conjecture}
\newtheorem{notation}[theorem]{Notation}

\newcommand{\R}{\mathbb R}
\newcommand{\Z}{\mathbb Z}
\newcommand{\C}{\mathbb C}
\newcommand{\K}{\mathbb K}
\newcommand{\N}{\mathbb N}
\newcommand{\bbQ}{\mathbb Q}
\newcommand{\Hyp}{\mathbb H}
\newcommand{\bbV}{\mathbb V}
\newcommand{\bbLambda}{\reflectbox{\raisebox{\depth}{\scalebox{1}[-1]{$\mathbb V$}}}}
\newcommand{\gl}{\mathfrak{gl}}
\newcommand{\sldue}{\mathfrak{sl}_2}
\newcommand{\so}{\mathfrak{so}}
\newcommand{\g}{\mathfrak g}
\newcommand{\X}{\mathfrak{X}}
\newcommand{\Int}{\mathrm{Int}}
\newcommand{\Ad}{\mathrm{Ad}}
\newcommand{\ad}{\mathrm{ad}}
\newcommand{\rank}{\mathrm{rank}}
\newcommand{\GL}{\mathrm{GL}}
\newcommand{\FL}{\mathrm{FL}}
\newcommand{\PSL}{\mathrm{PSL}}
\newcommand{\PSO}{\mathrm{PSO}}
\newcommand{\PGL}{\mathrm{PGL}}
\newcommand{\Bl}{\mathrm{Bl}}
\newcommand{\SL}{\mathrm{SL}}
\newcommand{\SO}{\mathrm{SO}}
\newcommand{\Homeo}{\mathrm{Homeo}}
\newcommand{\Diff}{\mathrm{Diff}}
\newcommand{\Isom}{\mathrm{Isom}}
\newcommand{\Aff}{\mathrm{Aff}}
\newcommand{\Ext}{\mathrm{Ext}}
\newcommand{\Tor}{\mathrm{Tor}}
\newcommand{\Aut}{\mathrm{Aut}}
\newcommand{\Hom}{\mathrm{Hom}}
\newcommand{\End}{\mathrm{End}}
\newcommand{\Rep}{\mathrm{Rep}}
\newcommand{\Vol}{\mathrm{Vol}}
\newcommand{\Lie}{\mathrm{Lie}}
\newcommand{\Span}{\mathrm{Span}}
\newcommand{\Spin}{\mathrm{Spin}}
\newcommand{\Stab}{\mathrm{Stab}}
\newcommand{\Skew}{\mathrm{Skew}}
\newcommand{\Ann}{\mathrm{Ann}}
\newcommand{\Vect}{\mathrm{Vect}}
\newcommand{\Funct}{\mathrm{Funct}}
\newcommand{\Pf}{\mathrm{Pf}}
\newcommand{\id}{\mathrm{id}}
\newcommand{\ed}{\mathrm d}
\newcommand{\der}{\frac{\mathrm{d}}{\mathrm{d}t}}
\newcommand{\CW}{\mathrm{CW}}
\newcommand{\tr}{\mathrm{tr}}
\newcommand{\sign}{\mathrm{sign}}
\newcommand{\im}{\mathrm{im} \;}
\newcommand{\dR}{\mathrm{dR}}
\newcommand{\Rib}{\mathrm{Rib}}
\newcommand{\Ob}{\mathrm{Ob}}
\newcommand{\col}{\mathrm{col}}
\newcommand{\op}{\mathrm{op}}
\newcommand{\Zed}{\mathrm Z}
\newcommand{\A}{\mathrm A}
\newcommand{\En}{\mathrm N}
\newcommand{\Vsp}{\mathrm V}
\newcommand{\Set}{\mathrm{Set}}
\newcommand{\pro}{\mathrm{pro}}
\newcommand{\Man}{\mathrm{Man}}
\newcommand{\ind}{\mathrm{ind}}
\newcommand{\UC}{\mathrm{UC}}
\newcommand{\WRT}{\mathrm{WRT}}
\newcommand{\CGP}{\mathrm{CGP}}
\newcommand{\Cinf}{\mathcal C^{\infty}}
\newcommand{\Cat}{\mathcal C}
\newcommand{\Alg}{\mathcal A}
\newcommand{\Q}{\mathcal Q}
\newcommand{\B}{\mathcal B}
\newcommand{\D}{\mathcal D}
\newcommand{\F}{\mathcal F}
\newcommand{\G}{\mathcal G}
\newcommand{\V}{\mathcal V}
\newcommand{\curlyLambda}{\reflectbox{\raisebox{\depth}{\scalebox{1}[-1]{$\mathcal V$}}}}
\newcommand{\M}{\mathcal M}
\newcommand{\Nat}{\mathcal N}
\newcommand{\Esse}{\mathcal S}
\newcommand{\Acca}{\mathcal H}
\newcommand{\Ell}{\mathcal L}
\newcommand{\T}{\mathcal T}
\newcommand{\Bord}{\mathcal{B \! o \! r \mkern-4mu d}}
\newcommand{\mind}{\mathbf i}
\newcommand\restr[3][0]{{
  \raisebox{-#1pt}{$\left.
  \raisebox{#1pt}{$#2$}
  \vphantom{\big|}
  \right|_{#3}$}
  }}
\newcommand{\ddx}{\curvearrowright}
\newcommand{\dsx}{\curvearrowleft}
\newcommand{\bdx}{\raisebox{\depth}{\scalebox{1}[-1]{$\curvearrowright$}}}
\newcommand{\bsx}{\raisebox{\depth}{\scalebox{1}[-1]{$\curvearrowleft$}}}
\newcommand{\one}{\mathbbm 1}
\newcommand{\two}{$\langle 2 \rangle$}
\newcommand\numberthis{\addtocounter{equation}{1}\tag{\theequation}}
\newcommand{\colim}{\operatornamewithlimits{colim}}
\newcommand{\sqtensor}{\mbox{\footnotesize$\; \boxtimes \;$\normalsize}}
\newcommand{\psqtensor}{\mbox{\tiny$\; \boxtimes \;$\normalsize}}
\renewcommand{\leq}{\leqslant} 
\renewcommand{\geq}{\geqslant}
\renewcommand{\P}{\mathbb P}
\renewcommand{\labelenumi}{\textnormal{(}\textit{\roman{enumi}}\textnormal{)}}
\newcommand\minore[1]{{
  \overset{\scriptscriptstyle <}{#1}
  }}
\newcommand\maggiore[1]{{
  \overset{\scriptscriptstyle >}{#1}
  }}

\title[Quantum Invariants from Non-Semisimple Categories]{Quantum Invariants of 3-Manifolds Arising from Non-Semisimple Categories}
\author[M. De Renzi]{Marco De Renzi}
\address{Université Paris Diderot -- Paris 7, Paris, France} 
\email{marco.de-renzi@imj-prg.fr}

\begin{abstract}
 This survey covers some of the results contained in the papers by Costantino, Geer and Patureau \cite{cgp1} 
 and by Blanchet, Costantino, Geer and Patureau \cite{bcgp}. In the first one the authors construct two families of
 Reshetikhin-Turaev-type invariants of 3-manifolds, $\En_r$ and $\En^0_r$, using non-semisimple categories of representations of a 
 quantum version of $\sldue$ at a $2r$-th root of unity with $r \geq 2$. The secondary invariants $\En^0_r$ conjecturally 
 extend the original Reshetikhin-Turaev quantum $\sldue$ invariants.
 The authors also provide a machinery to produce invariants out of more general ribbon categories which can lack the semisimplicity 
 condition. In the second paper a renormalized version of $\En_r$ for 
 $r \not\equiv 0 \; (\mathrm{mod} \; 4)$ is extended to a TQFT, and connections with classical invariants such as the 
 Alexander polynomial and the Reidemeister torsion are found. In particular, it is shown that the use of richer categories 
 pays off, as these non-semisimple invariants are strictly finer than the original semisimple ones: 
 indeed they can be used to recover the classification of lens spaces, which Reshetikhin-Turaev invariants 
 could not always distinguish.
\end{abstract}

\maketitle

\section{Modular categories}  

A \textit{(strict) ribbon category} $\Cat$ is a (strict) monoidal category equipped with a braiding $c$, a twist 
$\vartheta$ and a compatible duality $(*,b,d)$. We will tacitly assume that all the ribbon categories we consider are strict. 
The \textit{category $\Rib_{\Cat}$ of ribbon graphs over $\Cat$} is the ribbon category whose objects are finite sequences 
$(V_1,\varepsilon_1), \ldots, (V_k,\varepsilon_k)$ where $V_i \in \Ob(\Cat)$ and $\varepsilon_i = \pm 1$ and whose morphisms
are isotopy classes of $\Cat$-colored ribbon graphs which are compatible with sources and targets. 

\begin{figure}[btph]\label{elementary}
\centering
\includegraphics{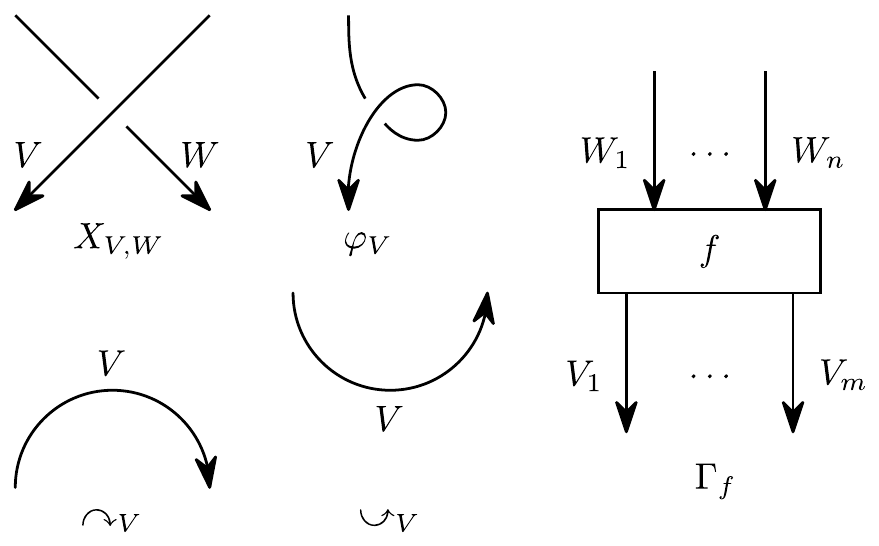}
\caption{Elementary $\Cat$-colored ribbon graphs}\label{cap_elementary}
\end{figure}

\begin{theorem}
 If $\Cat$ is a ribbon category then there exists a unique (strict) monoidal functor $F :\Rib_{\Cat} \rightarrow \Cat$ 
 such that (see Figure \ref{elementary}):
 \begin{enumerate}
  \item $F(V,+1) = V$ and $F(V,-1) = V^*$;
  \item $F(X_{V,W}) = c_{V,W}$, $F(\varphi_V) = \vartheta_V$, $F(\bdx_V) = b_V$ and $F(\ddx_V) = d_V$;
  \item $F(\Gamma_f) = f$.
 \end{enumerate}
\end{theorem}

The functor $F$ is the \textit{Reshetikhin-Turaev functor associated with $\Cat$}. 

\begin{remark}
 Every Reshetikhin-Turaev functor $F$ yields an invariant of framed oriented links colored with objects of $\Cat$.
\end{remark}

A \textit{ribbon Ab-category} is a ribbon category $\Cat$ whose sets of morphisms admit abelian group structures which make the 
composition and the tensor product of morphisms into $\Z$-bilinear maps. 
Then $K := \End_{\Cat}(\one)$ becomes a commutative ring called the \textit{ground ring of $\Cat$} and all sets of morphisms are 
naturally endowed with $K$-module structures (the scalar multiplication being given by tensor products with elements of $K$ on the 
left). 

An object $V \in \Ob(\Cat)$ is \textit{simple} if $\End(V) \simeq K$.

\begin{remark}
 We will always suppose that all the ribbon Ab-categories we consider have a field $\K$ for ground ring.
\end{remark}

A \textit{semisimple category} is a ribbon Ab-category $\Cat$ together with a distinguished set of simple objects 
$\Gamma(\Cat) := \{ V_i \}_{i \in I}$ such that:
\begin{enumerate}
 \item there exists $0 \in I$ such that $V_0 = \one$;
 \item there exists an involution $i \mapsto i^*$ of $I$ such that $V_{i^*} \simeq V_i^*$;
 \item for all $V \in \Ob(\Cat)$ there exist $i_1, \ldots, i_n \in I$ and maps $\alpha_j : V_{i_j} \rightarrow V$,
  $\beta_j : V \rightarrow V_{i_j}$ such that $\id_V = \sum_{j=1}^n \alpha_j \beta_j$ (we say that the set $\Gamma(\Cat)$
  \textit{dominates} $\Cat$);
 \item for any distinct $i, j \in I$ we have $\Hom_{\Cat}(V_i,V_j) = 0$.
\end{enumerate}

In a semisimple category we have the following results:

\begin{lemma}
 For all $V,W \in \Ob(\Cat)$:
 \begin{enumerate}
  \item $\Hom_{\Cat}(V,W)$ is a finite-dimensional $\K$-vector space.
  \item $\Hom_{\Cat}(V_i,V) = 0$ for all but a finite number of $i \in I$.
  \item $\Hom_{\Cat}(V,W) \simeq \bigoplus_{i \in I} \Hom_{\Cat}(V,V_i) \otimes_{\K} \Hom_{\Cat}(V_i,W)$ where the inverse 
   isomorphism is given by $f \otimes g \mapsto g \circ f$ on direct summands.
  \item the $\K$-bilinear pairing $\Hom_{\Cat}(V,W) \otimes_{\K} \Hom_{\Cat}(W,V) \rightarrow \K$ given by 
   $f \otimes g \mapsto \tr_{\Cat}(g \circ f)$ is non-degenerate.
 \end{enumerate}
\end{lemma}

\begin{corollary}
 The quantum dimension of simple objects is non-zero.
\end{corollary}

\begin{remark}
 Let $(f_i)_1, \ldots, (f_i)_{n_i}$ be a basis for the finite dimensional $\K$-vector space $\Hom_{\Cat}(V,V_i)$ and let 
 $(g_i)^1, \ldots (g_i)^{n_i}$ denote the dual basis of $\Hom_{\Cat}(V_i,V)$ defined by
 \[
  \dim_{\Cat}(V_i) (f_i)_h \circ (g_i)^k = \delta_h^k \cdot \id_{V_i},
 \]
 which exists thanks to (\textit{iv}) of the previous Proposition.
 Then we can write
 \[
  \id_V = \sum_{i \in I} \sum_{h,k = 1}^{n_i} (\lambda_i)^k_h \left[ (g_i)^k \circ (f_i)_h \right].
 \]
 But now we have
 \begin{align*}
  \delta_h^k \cdot \id_{V_i} &= \dim_{\Cat}(V_i) \left[ (f_i)_h \circ (g_i)^k \right] = 
  \dim_{\Cat}(V_i) \left[ (f_i)_h \circ \id_V \circ (g_i)^k \right] \\
  &= \dim_{\Cat}(V_i)^{-1} (\lambda_i)^k_h \cdot \id_{V_i}.
 \end{align*}
 Therefore
 \begin{equation}\label{fusion}
  \id_V = \sum_{i \in I} \sum_{j = 1}^{n_i} \dim_{\Cat}(V_i) \cdot (g_i)^j \circ (f_i)_j.
 \end{equation}
 Equation \ref{fusion} is called the \textit{fusion formula}.
\end{remark}

A \textit{premodular category} is a semisimple category $(\Cat,\Gamma(\Cat))$ such that $\Gamma(\Cat)$ is finite.

\begin{figure}[hbtp]\label{hopf}
\centering
\includegraphics{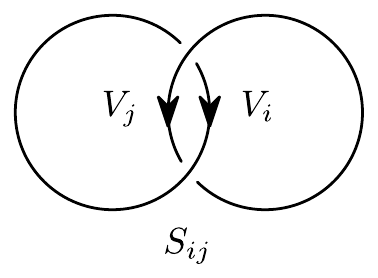}
\caption{Postive Hopf link colored with $V_i$ and $V_j$}\label{cap_hopf}
\end{figure}

A \textit{modular category} is a premodular category $(\Cat,\Gamma(\Cat)=\{ V_i \}_{i \in I})$ such that the matrix
$S = (F(S_{ij}))_{i,j \in I}$ with $S_{ij}$ given by Figure \ref{hopf} is invertible.

\section{The Reshetikhin-Turaev invariants}

The construction of Reshetikhin and Turaev associates with every premodular category $\Cat$ an invariant $\tau_{\Cat}$
of 3-manifolds (which will always be assumed to be closed and oriented)
provided $\Cat$ satisfies some non-degeneracy condition. Let us outline the general procedure in this context: 
let $\Cat$ be a premodular category, let $F : \Rib_{\Cat} \rightarrow \Cat$ be the associated Reshetikhin-Turaev functor and let
$\Omega$ be the associated \textit{Kirby color}
\[
 \Omega := \sum_{W \in \Gamma(\Cat)} \dim_{\Cat}(W) \cdot W.
\]

It is known that every 3-manifold $M^3$ can be obtained by surgery along some framed link $L$ inside $S^3$ (we write $S^3(L)$ for the 
result of this operation)
and that two framed
links yield the same 3-manifold if and only if they can be related by a finite sequence of Kirby moves. 
Therefore,
in order to find an invariant of 3-manifolds, we can look for an invariant of framed links which remains unchanged under
Kirby moves. For example let $L \subset S^3$ be a framed link giving a surgery presentation for $M^3$ and let
$\vec{L}(\Omega)$ denote the $\Cat$-colored ribbon graph obtained by assigning to each component an arbitrary orientation 
and the Kirby color $\Omega$. Then by evaluating $F$ on $\vec{L}(\Omega)$ we get a number in $\K$ and, thanks to the 
closure (up to isomorphism) of $\Gamma(\Cat)$ under duality, we can prove that $F(\vec{L}(\Omega))$ is actually
independent of the chosen orientation for $L$. Therefore we have a number $F(L(\Omega)) \in \K$ which depends only on the link 
$L$ giving a surgery presentation for $M$. Let us see its behaviour under Kirby moves.

\begin{proposition}[Slide]
 Let $(\Cat,\Gamma(\Cat))$ be a pre-modular category and let $T$ be a $\Cat$-colored ribbon graph. If 
 $T'$ is a $\Cat$-colored ribbon graph obtained from $T$ by performing a slide of an arc $e \subset T$ over a circle 
 component $K \subset T$ colored by $\Omega$, then $F(T') = F(T)$.
\end{proposition}

This result crucially relies on the semisimplicity of $\Cat$, which enables us to establish the fusion formula \ref{fusion},
and on the finiteness of $\Gamma(\Cat)$, which enables us to define Kirby colors.

Now let us turn our attention towards blow-ups and blow-downs. 
Let $\Delta_{\pm}$ denote the image under $F$ of a $\pm 1$-framed unknot colored by $\Omega$. 
If $L' \subset S^3$ is a link obtained from $L \subset S^3$ by a $\pm 1$-framed blow-up then 
$F(L'(\Omega)) = \Delta_{\pm} \cdot F(L(\Omega))$. At the same time we have 
that the positive and negative signatures of the linking matrices of $L'$ and $L$ satisfy
\[
 \sigma_{\pm}(L') = \sigma_{\pm} (L) + 1, \quad \sigma_{\mp}(L') = \sigma_{\mp} (L).
\]
Therefore we are tempted to consider the ratio
\[
 \frac{F(L(\Omega))}{\Delta_+^{\sigma_+(L)} \cdot \Delta_-^{\sigma_-(L)}},
\]
which is invariant under all Kirby moves. In order to be able to do so we must require from the premodular category $\Cat$ the 
following:

\medskip

\noindent{\textbf{Condition 1.} $\Delta_+ \cdot \Delta_- \neq 0$.}

\medskip

Therefore, let $\Cat$ be a premodular category satisfying Condition 1 and let $(M,T)$ be a pair consisting of a 
3-manifold $M$ and a closed $\Cat$-colored ribbon graph
$T \subset M$. If $L \subset S^3$ is any framed link yielding a surgery presentation for $M$ and $\Gamma_T$ is a $\Cat$-colored 
ribbon graph in $S^3 \smallsetminus L$ representing $T$ then the \textit{Reshetikhin-Turaev invariant associated with $\Cat$} is
\[
 \tau_{\Cat}(M,T) := \frac{F(L(\Omega) \cup \Gamma_T)}{\Delta_+^{\sigma_+(L)} \cdot \Delta_-^{\sigma_-(L)}}.
\]

\begin{remark}
 The actual Reshetikhin-Turaev invariant is given by the renormalization $\D^{- b_1(M) -1} \tau_{\Cat}(M)$ 
 where $b_1(M)$ is the first Betti number of $M$ and $\D$ is an element of $\K$ satisfying $\D^2 = F(u(\Omega))$ 
 with $u(\Omega)$ the $\Omega$-colored 0-framed unknot. Note that such a $\D$  may not exist and we may 
 have to manually adjoin it (compare with \cite{turaev-quantum}).
\end{remark}

\begin{remark}
 The non-degeneracy condition is automatically satisfied by any modular category.
\end{remark}

The most famous example of this construction, which yields the original invariants defined by Reshetikhin and Turaev, 
is obtained by considering a representation category 
of a quantum version of $\sldue$ at a root of unity. Let us recall the construction: fix an integer $r \geq 2$, 
set $q := e^{\frac{\pi}{r} i}$ and consider the quantum group $U_q(\sldue)$ generated (as a unital $\mathbb C$-algebra)
by $E,F,K,K^{-1}$ with relations
\begin{gather*}
 KK^{-1} = K^{-1}K = 1, \quad
 KEK^{-1} = q^2 E, \quad KFK^{-1} = q^{-2} F, \\
 [E,F] = \frac{K - K^{-1}}{q - q^{-1}}, \quad E^r=F^r=0
\end{gather*}
and with comultiplication, counit and antipode given by
\begin{gather*}
 \Delta(E) = E \otimes K + 1 \otimes E, \quad
 \varepsilon(E) = 0, \quad S(E) = -EK^{-1}, \\
 \Delta(F) = F \otimes 1 + K^{-1} \otimes F, \quad
 \varepsilon(F) = 0, \quad S(F) = -KF, \\
 \Delta(K^{\pm 1}) = K^{\pm 1} \otimes K^{\pm 1}, \quad
 \varepsilon(K^{\pm 1}) = 1, \quad S(K^{\pm 1}) = K^{\mp 1}.
\end{gather*}
A representation of $U_q(\sldue)$ is a \textit{weight representation}, or a weight $U_q(\sldue)$-module, 
if it splits as a direct sum of eigenspaces for the action of $K$.
The Hopf algebra structure on $U_q(\sldue)$ endows the category $U_q(\sldue)$-mod of finite-dimensional weight 
representations of the quantum group $U_q(\sldue)$ with a natural monoidal structure and a compatible duality. 
Now let $\bar{U}_q(\sldue)$ denote the quantum group obtained from $U_q(\sldue)$ by adding 
the relation $K^{2r} = 1$. This condition forces all weights (eigenvalues for the action of $K$) 
to be integer powers of $q$ for all representations of $\bar{U}_q(\sldue)$. Therefore we can consider
the operator $q^{H \otimes H / 2}$ defined on $V \otimes W$ for all weight $\bar{U}_q(\sldue)$-modules $V$ and $W$
by the following rule:
\[
 q^{H \otimes H / 2}(v \otimes w) = q^{m n/2} v \otimes w 
\]
if $Kv = q^m v$ and $Kw = q^n w$, where $q^{m n/2}$ stands for $e^{\frac{mn \pi}{2r} i}$.
Set $\{ m \} := q^m - q^{- m}$ for all $m \in \Z$ and define
\[
 [n] := \frac{\{ n \}}{\{ 1 \}}, \quad [n]! := [n] [n-1] \cdots [1].
\]
for all $n \in \N$. Consider the operator $R$ defined on $V \otimes W$ for all weight $\bar{U}_q(\sldue)$-modules $V$ and $W$ as
\[
 q^{H \otimes H/2} \sum_{n=0}^{r-1} \frac{q^{n(n-1)/2}}{[n]!} \{ 1 \}^n E^n \otimes F^n.
\]
Finally consider the operator $q^{-H^2 / 2}$ determined on each weight $\bar{U}_q(\sldue)$-module $V$ by the rule
\[
 q^{-H^2 / 2}(v) = q^{-n^2/2} v \mbox{ if } Kv = q^n v,
\]
define the operator $u$ as 
\[
 q^{-H^2/2} \sum_{n=0}^{r-1} \frac{q^{3n(n-1)/2}}{[n]!} \{ -1 \}^n F^n K^{-n} E^n 
\]
and set $v := K^{r-1} u$.
Then the category $\bar{U}_q(\sldue)$-mod of finite dimensional weight representations of $\bar{U}_q(\sldue)$ 
can be made into a ribbon Ab-category by considering the compatible braidings and 
twists given by
\[
 c_{V,W} = \tau \circ R : V \otimes W \rightarrow W \otimes V, \quad
 \vartheta_V = v^{-1} : V \rightarrow V, 
\]
where $\tau$ is the $\K$-linear map switching the two factors of the tensor product. 
Moreover $\bar{U}_q(\sldue)$-mod is quasi-dominated by a finite number of simple modules, and thus it
can be made into a modular category by quotienting negligible morphisms. The invariant we obtain is denoted $\tau_r$.

\section{Relative G-premodular categories}

To motivate the construction of non-semisimple invariants, let us consider the following different quantization of $\sldue$:
let $U^H_q(\sldue)$ denote the quantum group obtained by adding to $U_q(\sldue)$ the additional generator $H$ satisfying
the following relations:
\begin{gather*}
 HK = KH, \quad [H,E]=2E, \quad [H,F]=-2F \\
 \Delta(H) = H \otimes 1 + 1 \otimes H, \quad
 \varepsilon(H) = 0, \quad S(H) = -H.
\end{gather*}

\begin{remark}
 The new generator $H$ should be thought of as a logarithm of $K$ and, even though we will not require the relation to hold true
 at the quantum group level, we will restrict ourselves to representations where it is satisfied.
\end{remark}

The category $U^H_q(\sldue)$-mod of finite-dimensional weight representations of the quantum group $U^H_q(\sldue)$ 
where $K$ acts like the operator $q^H$
can be made into a ribbon Ab-category by means of the same R-matrix and ribbon element used for 
$\bar{U}_q(\sldue)$-mod.

\begin{remark}
 The introduction of $H$ is necessary in order to make sense of the formulas defining the operators $R$ and $u$ because the 
 absence of the relation $K^{2r} = 1$ makes room for weights which are not integer powers of $q$. 
 The operator $q^{H \otimes H / 2}$ is then given by the rule
 \[
  q^{H \otimes H / 2}(v \otimes w) = q^{\lambda \mu/2} v \otimes w 
 \]
 if $Hv = \lambda v$ and $H w = \mu w$, where $q^{\alpha}$ stands for $e^{\frac{\alpha \pi}{r} i}$ for all $\alpha \in \mathbb C$.
 The definition of $q^{-H^2/2}$ is analogous.
\end{remark}

What is different in $U^H_q(\sldue)$-mod is that simple objects are not in a finite number: 
indeed for any $\alpha \in (\mathbb C \smallsetminus \Z) \cup r \cdot \Z$ 
the $r$-dimensional module $V_{\alpha}$ generated by the highest weight vector $v_0^{\alpha}$ satisfying $E v_0^{\alpha} = 0$ and 
$H v_0^{\alpha} = (\alpha + r -1) v_0^{\alpha}$ is simple and projective, and is called a \textit{typical module} (see 
\cite{cgp2} for details). 
If we could put this richer category into the Reshetikhin-Turaev machinery we would perhaps find more refined 3-manifold invariants. 
It is indeed the case, but we need to face (among other things) the following obstructions:
\begin{enumerate}
 \item every typical module has zero quantum dimension;
 \item we cannot quotient negligible morphisms as this would kill all typical modules, and thus we are forced to work
  with a non-semisimple category;
 \item typical modules are pairwise non-isomorphic, and therefore we have to deal with infinitely many isomorphism classes of 
  simple objects.
\end{enumerate}
Let us see how we can work around these obstacles. The idea is to generalize the Reshetikhin-Turaev construction
to more general ribbon Ab-categories which have the previous set of obstructions.

\subsection*{Facing obstruction (\textit{i}): Modified quantum dimension}

To begin with let us take care of the vanishing quantum dimension problem.
The strategy will be to use categories $\Cat$ which admit a modified dimension which does not vanish. In order to do so
we need an \textit{ambidextrous pair $(\A,\ed)$}, that is the given of a set of simple objects $\A \subset \Ob(\Cat)$ and a map 
$\ed : \A \rightarrow \K^*$ with the following property:
if $T$ is an \textit{$\A$-graph}, i.e. a closed $\Cat$-colored ribbon graph admitting at least one color in $\A$, if $e \subset T$ is
an arc colored by $V \in \A$ and if $T_e$ denotes the element of $\End_{\Rib_{\Cat}}((V,+))$ obtained by cutting open $T$ at $e$, 
then 
\[
 F'(T) := \ed(V) \cdot \langle T_e \rangle
\]
is independent of the chosen $\A$-colored arc $e$ (here $\langle T_e \rangle$ denotes the unique element of $\K$ such that
$F(T_e) = \langle T_e \rangle \cdot \id_V$). 

\begin{definition}
 A ribbon Ab-category $\Cat$ admitting an ambidextrous pair $(\A,\ed)$ is said to have \textit{modified dimension $\ed$}.
 $F'$ is the \textit{modified $\A$-graph invariant associated with $(\A,\ed)$}.
\end{definition}

\begin{example}
 In the category $U^H_q(\sldue)$-mod considered before we have indeed an ambidextrous pair. It is obtained by taking $\A$ to be 
 the set of typical modules and by defining
 \[
  \ed(V_{\alpha}) = (-1)^{r-1} \prod_{j=1}^{r-1} \frac{q^j - q^{-j}}{q^{\alpha + r - j} - q^{-\alpha - r + j}} 
                = (-1)^{r-1} \cdot r \cdot \frac{\sin(\alpha \pi / r)}{\sin(\alpha \pi)}.
 \]
\end{example}

\subsection*{Facing obstruction (\textit{ii}): G-grading relative to X}

Moving on to the subject of semisimplicity, we will ask our categories to have a 
distinguished family of semisimple full subcategories arranged into a grading, i.e. indicized by an abelian group $G$ in such a way 
that the tensor product respects the group operation. The aim of course is to work as much as possible in the graded semisimple part
of the category and to leave aside the non-semisimple part.

\begin{definition}
 Let $\Cat$ be a ribbon Ab-category. A full subcategory $\Cat'$ of $\Cat$ is said to be \textit{semisimple inside $\Cat$} if
 it is dominated by a set $\Gamma(\Cat')$ of simple objects of $\Cat'$ such that for any distinct $V,W \in \Gamma(\Cat')$ we have 
 $\Hom_{\Cat}(V,W) = 0$.
\end{definition}

\begin{remark}
 We do not ask of $\Gamma(\Cat')$ to contain $\one$ nor to be closed under duality up to isomorphism. In particular it may
 very well happen that the quantum dimension of a simple object of $\Cat'$ is zero.
\end{remark}

\begin{definition}
 We will say that a subset $X$ of an abelian group $G$ is \textit{small} if $G$ cannot be covered by any finite union of
 translated copies of $X$, i.e. if there exists no choice of $g_1, \ldots, g_k \in G$ such that
 \[
  G \subset \bigcup_{i=1}^k (g_i + X).
 \]
 Let $G$ be an abelian group and let $X \subset G$ be a small subset. 
 A family of full subcategories $\{ \Cat_g \}_{g \in G}$ of a ribbon Ab-category category 
 $\Cat$ gives a \textit{$G$-grading relative to $X$} for 
 $\Cat$ if:
 \begin{enumerate}
  \item $\Cat_g$ is semisimple inside $\Cat$ for all $g \in G \smallsetminus X$;
  \item $V \in \Ob(\Cat_g), \; V' \in \Ob(\Cat_{g'}) \Rightarrow V \otimes V' \in \Ob(\Cat_{g + g'})$;
  \item $V \in \Ob(\Cat_g) \Rightarrow V^* \in \Ob(\Cat_{- g})$;
  \item $V \in \Ob(\Cat_g), \; V' \in \Ob(\Cat_{g'}), \; g \neq g' \Rightarrow \Hom_{\Cat}(V,V') = 0$.
 \end{enumerate}
\end{definition}

The elements of $g$ which are not contained in $X$ are called \textit{generic} and a subcategory $\Cat_g$ indicized by a
generic $g$ is called a \textit{generic subcategory}.
A category $\Cat$ with a $G$-grading relative to $X$ will be called a \textit{$G$-graded category} for the sake of brevity.

\begin{example}
 In the category $U^H_q(\sldue)$-mod considered before we have a relative $G$-grading too. Indeed we can take 
 $G = \mathbb C / 2 \Z$, $X = \Z / 2 \Z$ and set $\Cat_{\bar{\alpha}}$ 
 equal to the full subcategory of modules whose weights are all
 congruent to $\alpha$ modulo 2. Then every $\Cat_{\bar{\alpha}}$ with $\alpha$ not integer is semisimple inside 
 $\bar{U}_q(\sldue)$-mod, being dominated by the typical modules it contains.
\end{example}

\subsection*{Facing obstruction (\textit{iii}): Periodicity group}

Finally, for the finiteness issue, we will proceed as follows: for a $G$-graded category $\Cat$ 
we will ask the sets of isomorphism classes of simple objects of all generic subcategories
to be finitely partitioned in a way we can control. 

\begin{definition}
 A set $C \subset \Ob(\Cat)$ of objects of a ribbon Ab-category is a \textit{commutative family} 
 if the braiding and the twist are trivial
 on $C$, i.e. if we have $c_{W,V} \circ c_{V,W} = \id_{V \otimes W}$ and $\vartheta_V = \id_V$ for all $V,W \in C$.
\end{definition}

\begin{definition}
 Let $\Zed$ be an abelian group and $\Cat$ be a ribbon Ab-category. A \textit{realization} of $\Zed$ in $\Cat$ is a commutative
 family $\{ \varepsilon^t \}_{t \in \Zed}$ satisfying
 \[
  \varepsilon^0 = \one, \quad \varepsilon^t \otimes \varepsilon^s = \varepsilon^{t+s}, \quad
  \dim_{\Cat}(\varepsilon^t) = 1 \quad \forall
  \; t,s \in \Zed.
 \]
\end{definition}

Any free realization of $\Zed$ gives isomorphisms between the $\K$-vector spaces of morphisms $\Hom_{\Cat}(V,W)$ and 
$\Hom_{\Cat}(V \otimes \varepsilon^t,W \otimes \varepsilon^t)$ for all choices of $V,W \in \Ob(\Cat)$ and $t \in \Zed$.
Indeed the inverse of the map $f \mapsto f \otimes \id_{\varepsilon^t}$ is simply given by 
$g \mapsto g \otimes \id_{\varepsilon^{-t}}$.
Therefore if $V$ is simple then $V \otimes \varepsilon^t$ is simple too for all $t \in \Zed$. 
Thus any realization of $\Zed$ induces an action of $\Zed$ on (isomorphism classes of) objects of $\Cat$ given by the tensor 
product on the right with $\varepsilon^t$. Such a realization is \textit{free} if this action is free.

\begin{definition}
 An abelian group $\Zed$ is the \textit{periodicity group} of the $G$-graded category $\Cat$ if there exists a free realization
 of $\Zed$ in $\Cat_0$ whose action on $\Gamma(\Cat_g)$ has a finite number of orbits for all $g \in G \smallsetminus X$.
\end{definition}

In this case there exists some finite set of representatives of $\Zed$-orbits $O(\Cat_g) \subset \Gamma(\Cat_g)$ 
for all generic $g$ such that
each simple module in $\Gamma(\Cat_g)$ is isomorphic to some tensor product $W \otimes \varepsilon^t$
for $W \in O(\Cat_g)$ and $t \in \Zed$.

\begin{example}
 Once again the category $U^H_q(\sldue)$-mod considered before gives us an instance of this structure. 
 Namely the periodicity group is $\Zed = \Z$ and its free realization in $\Cat_{\bar{0}}$ is given by a 1-dimensional 
 module $\varepsilon^t$ for every $t \in \Z$ which is spanned by the non-zero vector $v^t$ such that
 $Ev^t = Fv^t = 0$ and $Hv^t = 2rt v^t$.
\end{example}

The categories which will allow us to extend the Reshetikhin-Turaev construction to the non-semisimple
case admit all of the structures we just introduced.

\begin{figure}[btph]
\centering
\includegraphics{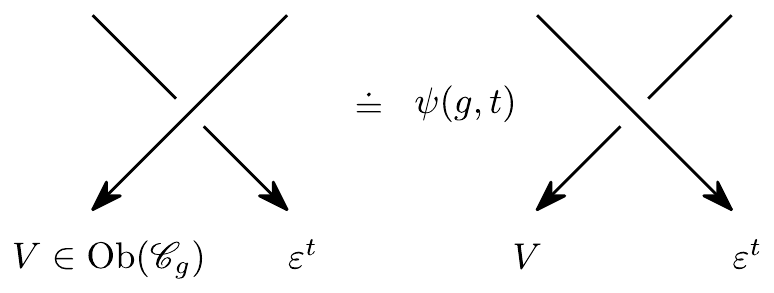}\label{skein}
\caption{Skein-type relation for $G$ and $\Zed$ (the $\doteq$ sign stands for equality under $F$)}
\end{figure}

\begin{definition}
 A \textit{relative $G$-premodular category} is $(\Cat,G \supset X,(\A,\ed),\Zed)$ where $\Cat$ is a $G$-graded category
 with modified dimension $\ed$ and periodicity group $\Zed$ satisfying the following compatibility conditions:
 \begin{enumerate}
  \item $\A \supset \Gamma(\Cat_g)$ for all $g \in G \smallsetminus X$;
  \item $c_{V,\varepsilon^t} = \psi(g,t) \cdot c_{\varepsilon^t,V}^{-1}$ for all 
   $V \in \Ob(\Cat_g), t \in \Zed$ and for some $\Z$-bilinear pairing $\psi : G \times \Zed \rightarrow \K^*$ (see Figure 
   \ref{skein}).
 \end{enumerate}
\end{definition}

\begin{example}
 $U^H_q(\sldue)$-mod is a relative $\mathbb C / \Z$-premodular category as it can be shown that a skein relation like the one 
 required in condition $($ii$\,)$ of the previous definition holds.
\end{example}

\section{Construction of non-semisimple invariants}

We are ready to sketch a construction analogous to the one of Reshetikhin and Turaev which associates with each 
relative $G$-premodular category $\Cat$ an invariant of 3-manifolds provided $\Cat$ satisfies some
non-degeneracy conditions.
The idea will be to use the modified invariant $F'$ as a basis for this construction exactly
as the functor $F$ was used as a basis for the standard case. Remember however that in order to compute $F'$ on a $\Cat$-colored
ribbon graph $T$ we will need to make sure that $T$ is actually an A-graph. 

Let us fix a relative $G$-premodular category $\Cat$. The first thing we did in the construction of Reshetikhin-Turaev 
invariants was to color a framed link giving a surgery presentation for a 3-manifold $M$ with the Kirby color $\Omega$ associated
with some premodular category. Now, in $\Cat$ we do not 
have the concept of a Kirby color, but we can define an infinite family of modified Kirby colors. 

Indeed if $g \in G$ is generic
then the formal sum 
\[
 \Omega_g := \sum_{W \in O(\Cat_g)} \ed(W) \cdot W
\]
is a \textit{modified Kirby color of degree $g$}.

\begin{remark}
 It can be easily proved using the properties of the periodicity group $\Zed$ that the modified dimension $\ed$ factorizes through
 a map defined on $\Zed$-orbits on all generic subcategories, i.e. we have $\ed(W \otimes \varepsilon^t) = \ed(W)$
 for all $W \in O(\Cat_g)$ and all $t \in \Zed$. In particular the coefficients in the formal sum $\Omega_g$ are independent of the 
 choice of the representatives of $\Zed$-orbits in $\Gamma(\Cat_g)$. Of course $W$ and $W \otimes \varepsilon^t$ are not 
 isomorphic if $t \neq 0$ but we will see that under certain circumstances this choice will not affect the value of $F'$.
\end{remark}

Since we defined an infinite family of modified Kirby colors it is not clear which one should be used to color the components
of a surgery link $L$ for $M^3$. The right choice is to let the coloring be determined by a cohomology class 
\[
 \omega \in 
H^1(M \smallsetminus T; G) \cong \Hom_{\Z}(H_1(M \smallsetminus T),G)
\]
which is compatible with the $\Cat$-coloring which is already present on $T$.

\begin{definition}
 Let $T$ be a $\Cat$-colored ribbon graph inside $M$ and $\omega$ be an element of $H^1(M \smallsetminus T; G)$. 
 For every arc $e \subset T$ let $\mu_e$ denote the homology class of a positive meridian around $e$. The triple $(M,T,\omega)$
 is \textit{compatible} if the color of $e$ is an object of $\Cat_{\langle \omega , \mu_e \rangle}$.
\end{definition}

We will now look for an invariant of compatible triples $(M,T,\omega)$, where two triples 
$(M_i,T_i,\omega_i)$ for $i = 1,2$ are considered to be equivalent if there exists an orientation preserving diffeomorphism
$f: M_1 \rightarrow M_2$ such that $f(T_1) = T_2$ as $\Cat$-colored ribbon graphs and $f^*(\omega_2) = \omega_1$.

\begin{remark}
 We will have to be more careful and to keep track of the (isotopy class of the) diffeomorphism induced by each Kirby move.
\end{remark}

The idea is to color each component $L_i$ of a surgery link $L$ with a modified Kirby color whose degree is 
determined by the evaluation $\langle \omega, \mu_i \rangle$, where $\mu_i$ denotes the homology class corresponding to a 
positive meridian of $L_i$. Thus, since modified Kirby colors are defined only for generic degrees, not all surgery 
presentations can be used to define the new invariant.

\begin{definition}
 A compatible triple $(M,T,\omega)$ admits a \textit{computable surgery presentation} 
 $L = L_1 \cup \ldots \cup L_m \subset S^3$ if one of the following holds:
 \begin{enumerate}
  \item $L \neq \varnothing$ and $\langle \omega, \mu_i \rangle$ is generic for all $i = 1, \ldots, m$;
  \item $L = \varnothing$ and $T$ is an $A$-graph.
 \end{enumerate}
\end{definition}

If $L$ is a link yielding a computable surgery presentation for a compatible triple $(M,T,\omega)$ and we denote by $L(\omega)$ the 
$\Cat$-colored link obtained by coloring each component $L_i$ of $L$ with $\Omega_{\langle \omega,\mu_i \rangle}$, then $F'$ can 
be evaluated on $L(\omega) \cup \Gamma_T$, where $\Gamma_T$ represents $T$ inside $S^3 \smallsetminus L$. 

\begin{remark}\label{computable}
 It can be shown that a sufficient condition for the existence of a computable surgery presentation for a compatible triple 
 $(M,T,\omega)$ is that the image of $\omega$ is not entirely contained in the critical set $X$ (when we regard $\omega$ as a 
 map from $H_1(M \smallsetminus T)$ to $G$).
\end{remark}

Let us see what happens when we perform Kirby moves.

\begin{definition}
 Let $T$ be a ribbon graph. A $G$-coloring for $T$ is a homology class $\varphi \in H_1(T;G)$. A $\Cat$-coloring
 on $T$ is \textit{$\varphi$-compatible} if for all arc $e \subset T$ the color of $e$ is an element of $\Cat_{g(e)}$ with
 \[
  \left[ \sum_{e \subset T} g(e) \cdot e \right] = \varphi.
 \]
\end{definition}

\begin{remark}
 Every cohomology class $\omega \in H^1(M \smallsetminus T; G)$ determines a $G$-coloring $\varphi_{\omega}$ for $L \cup \Gamma_T$.
\end{remark}

\begin{definition}
 Let $T \subset S^3$ be a $\varphi$-compatible $\Cat$-colored ribbon graph for some $\varphi \in H_1(T;G)$ and let 
 $K \subset S^3 \smallsetminus T$ be a framed knot. Then the \textit{$G$-linking number between $K$ and $\varphi$} is defined as
 \[
  \ell k^G(K,\varphi) := [\Sigma_{K}] \pitchfork \iota_*(\varphi)
 \]
 where 
 \[
  \pitchfork : H_2(S^3 \smallsetminus N(K),\partial N(K);G) \times H_1(S^3 \smallsetminus N(K);G) \rightarrow G
 \]
 is the intersection form with $G$ coefficients, $N(K)$ is a tubular neighborhood of $K$ disjoint from $T$,
 $\Sigma_K$ is a Seifert surface for a parallel copy of $K$ determined by the framing and 
 $\iota : T \hookrightarrow S^3 \smallsetminus N(K)$ denotes inclusion.
\end{definition}

\begin{remark}
 It can be shown that if $L \subset S^3$ is a computable surgery presentation for the compatible triple $(M,T,\omega)$ 
 and if $L'_i \subset S^3 \smallsetminus (L \cup \Gamma_T)$ is a parallel copy of a component of $L_i \subset L$ determined
 by the framing then 
 \[
  \ell k^G(L'_i,\varphi_{\omega}) = 0.
 \]
\end{remark}

\begin{figure}[hbtp]
\centering
\includegraphics{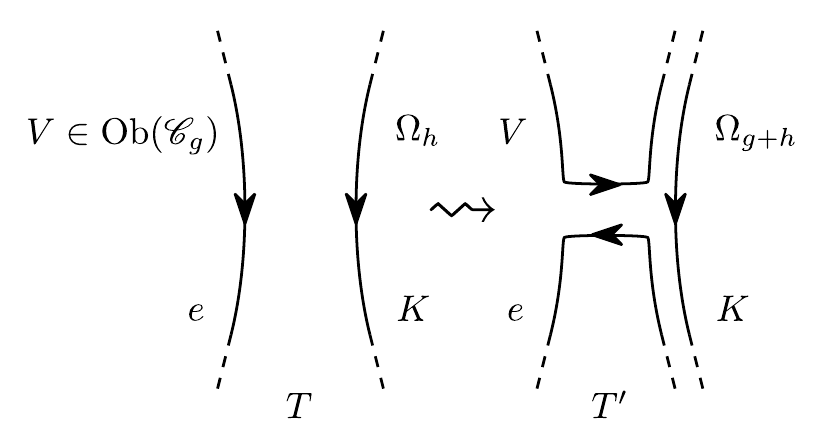}\label{slide_2}
\caption{Subtraction}
\end{figure}

\begin{proposition}[Slide]
 Let $T$ be a $\varphi$-compatible $\A$-graph, let $e \subset T$ be an arc colored by $V \in \Ob(\Cat_g)$ and let $K \subset T$ be a knot 
 component colored by $\Omega_h$ for some generic $h \in G$. Suppose that $g + h$ is generic too and that
 $\ell k^G(K',\varphi) = 0$ where $K'$ is a parallel copy of $K$ determined by the framing. If $T'$ is an 
 $\A$-graph obtained from $T$ by sliding $e$ along $K$ and switching the color of $K$ to 
 $\Omega_{g + h}$ (like in Figure \ref{slide_2}) then $F'(T') = F'(T)$.
\end{proposition}

To prove this proposition we need to establish a fusion formula (which can be done in the semisimple part of $\Cat$ exactly as
before) and to use the skein-type relation in the definition of $\Cat$ in order to handle closed components colored with 
$\varepsilon^t$ for $t \in \Zed$. The color of $K$ changes because $V \otimes W$ is an object of $\Cat_{g+h}$ for all 
$W \in O(\Cat_h)$.

For what concerns blow-ups and blow-downs, we cannot compute $F'$ directly on a detached $\pm 1$-framed unknot as such a component
should be colored with the modified Kirby color of degree 0 and it may very well 
happen that $0 \in X$ (which is the case in our previous example).

\begin{figure}[htbp]
\centering
\includegraphics{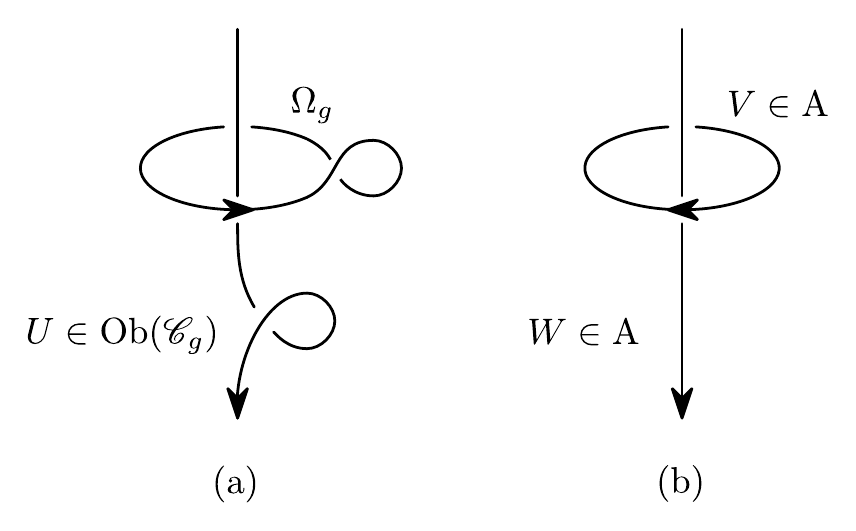}\label{blow-up_H-stab}
\caption{Blow-up of $+1$-framed meridian (a) and $H$-stabilization (b)}
\end{figure}

\begin{proposition}[Blow-up and blow-down]
 Let $T_+$ be the $\Cat$-colored ribbon graph given by Figure \ref{blow-up_H-stab}(a) with $g$ generic in $G$. 
 Then $\Delta_+ := \langle T_+ \rangle$ does not depend 
 on the generic $g$ nor on the object $U \in \Ob(\Cat_g)$. The same holds for the analogous graph $T_-$ (obtained by turning each 
 overcrossing of $T_+$ into an undercrossing) and for $\Delta_- := \langle T_- \rangle$.
\end{proposition}

\begin{remark}
 The operation of blowing up a positive meridian of an arc in a ribbon graph $T$ can replace the operation of blowing up
 an isolated unknotted componend provided $T$ is non-empty. 
 This is always the case for computable surgery presentations since there is always
 at least one arc colored in $\A$.
\end{remark}

Thus what we need in order to be able to define the invariant is once again to ask the condition $\Delta_+ \cdot \Delta_- \neq 0$.
However, this time we need also another non-degeneracy condition which allows us to perform an operation called 
$H$-stabilization which is needed in the proof of the invariance of our construction. 
Namely, let $H(V,W)$ denote the $\Cat$-colored ribbon graph given by Figure \ref{blow-up_H-stab}(b) for $V,W \in \A$. Then:

\medskip

\noindent{\textbf{Condition 2.} \begin{enumerate}
                                 \item \textit{$\Delta_+ \cdot \Delta_- \neq 0$}
                                 \item \textit{$\langle H(V,W) \rangle \neq 0$ for all $V,W \in \A$.}
                                \end{enumerate}
}

\medskip

Now we can state our result.

\begin{theorem}\label{teorema}
 Let $\Cat$ be a relative $G$-premodular category satisfying the non-degeneracy Condition 2. 
 Let $L$ be a framed link giving a computable surgery presentation for a compatible triple $(M,T,\omega)$ and let
 $\Gamma_T$ be a $\Cat$-colored ribbon graph inside $S^3 \smallsetminus L$ representing $T$. Then
 \[
  \En_{\Cat} (M,T,\omega) := \frac{F'(L(\omega) \cup \Gamma_T)}{\Delta_+^{\sigma_+(L)} \cdot \Delta_-^{\sigma_-(L)}}
 \]
 is a well-defined invariant of $(M,T,\omega)$.
\end{theorem}

\begin{remark}
 When $\Cat = U^H_q(\sldue)$-mod with $q = e^{\frac{\pi i}{r}}$ we write $\En_r$ instead of 
 $\En_{U^H_q(\sldue)\mathrm{-mod}}$.
\end{remark}

The subtlety in the proof of this result is the following: 
if we have two different computable surgery presentations $L$ and $L'$
it may happen that the sequence of Kirby moves connecting them passes through some non-computable presentation.
What we have to prove is that, up to passing to a different sequence of Kirby moves, we can make sure to get a
computable presentation at each intermediate step. 

This turns out to be true, but we have to allow for an operation, 
called $H$-stabilization, which modifies the triple $(M,T,\omega)$ and which is defined as follows:
let $e \subset T$ be an arc colored by $W \in \A$, let $\alpha$ be a positive 0-framed meridian of $e$ 
disjoint from $T$ and colored by $V \in \Gamma(\Cat_g)$ for some generic $g$ and let $D^2 \subset S^3$ be a disc  
intersecting $e$ once and satisfying $\partial D^2 = \alpha$. Now let $T_H$ denote the $\A$-graph $T \cup \alpha$ and let $\omega_H$ 
be the cohomology class coinciding with $\omega$ on homology classes contained in $M \smallsetminus (T \cup D^2)$ 
and satisfying $\langle \omega_H , \mu_{\alpha} \rangle = g$ where $\mu_{\alpha}$ is the homology class of a positive meridian of 
$\alpha$. 
Then the compatible triple $(M,T_H,\omega_H)$ is said to be obtained by \textit{$H$-stabilization of degree $g$} 
from $(M,T,\omega)$, and $\alpha$ is called the \textit{stabilizing meridian}. 
Now $(M,T_H,\omega_H)$ is not equivalent to $(M,T,\omega)$ but we have the equality
$\En_{\Cat} (M,T_H,\omega_H) = \langle H(V,W) \rangle \cdot \En_{\Cat} (M,T,\omega)$. 

Returning to the proof of the Theorem, we split the argument into three steps:
we begin by first proving the result in the case that $T$ itself is an $\A$-graph, that the initial and final surgery presentations
are the same and that the sequence of Kirby moves involves only isotopies of $\Gamma_T$ inside $S^3(L)$, i.e. slides
of arcs of $\Gamma_T$ over components of $L$ (we call this sequence of moves an \textit{isotopy inside $S^3(L)$}).
This case can be easily treated by performing a single $H$-stabilization on $(M,T,\omega)$ whose degree is 
sufficiently generic. Indeed if we slide a stabilizing meridian on some component $L_j$ of the computable link $L$ which is colored 
by $\Omega_{h_j}$ we change the color of $L_j$, 
provided the degree $g$ of the $H$-stabilization satisfies $g + h_j \in G \smallsetminus X$.
This would impose a condition on the choice of the degree, but there surely exists a $g \in G$ which satisfies it because 
$X$ is small. More generally, if $C \subset G$ denotes the finite set of (degrees of) colors appearing on $L$
during the sequence of slides, we can choose the degree $g$ of the $H$-stabilization in such a way that $(g + C) \cap X = 
\varnothing$. Thus, we can begin by sliding the stabilizing meridian $\alpha$ over all components of $L$, 
then we can follow the original sequence of Kirby moves and finally we can slide back $\alpha$ to its original position. 
What we will get is an equality of the form
\[
 \frac{F'(L(\omega) \cup \Gamma_T) \cdot \langle H(V,W) \rangle}{\Delta_+^{\sigma_+(L)} \cdot \Delta_-^{\sigma_-(L)}}
 = \; \frac{F'(L(\omega) \cup \Gamma'_T) \cdot \langle H(V,W) \rangle}{\Delta_+^{\sigma_+(L)} \cdot 
 \Delta_-^{\sigma_-(L)}},
\]
for some $V, W \in \A$, which proves the first step.

The second step consists in proving the Theorem when $T$ is an $\A$-graph. In this case, if
\[
 \begin{tikzcd}[column sep=15pt,row sep=-5pt]
  L^0 \cup \Gamma^0_T \arrow{r}{s_1} & \cdots \arrow{r}{s_k} & L^k \cup \Gamma^k_T \\
  \shortparallel & & \shortparallel \\
  L \cup \Gamma_T & & L' \cup \Gamma'_T
 \end{tikzcd}
\]
is our sequence of Kirby moves, we perform an $H$-stabilization for each 
Kirby move $s_h$ which makes some non-generic color appear. If $s_h$ is a non-admissible slide over some component $L^{h-1}_j$ 
we precede it by a slide of the corresponding stabilizing meridian $\alpha_h$ over $L^{h-1}_j$. If $s_h$ is a non-admissible 
blow-up around some arc we perform it on the 
corresponding stabilizing meridian $\alpha_h$ instead and then we slide the arc over the newly created component. 
All degrees can be chosen
so to adjust all colors, and the use of different stabilizations ensures the independence of the conditions. The tricky point is 
that a move $s_{\ell}$ which in the original sequence was a blow-down of a $\pm 1$-framed meridian of some arc or link component 
may now have become the blow-down of a component which is also linked to some of the stabilizing meridians we added. 
In this case, though, we can slide all these stabilizing meridians off, and this operation is an isotopy inside $S^3(L^{\ell - 1})$. 
Remark that the configuration we get at this point is not necessarily admissible, 
but problems can arise only for blow-downs of
meridians of arcs in $\Gamma_T^{\ell - 1}$. Thus in this case we can perform a new $H$-stabilization, 
slide the arc off and slide the new stabilizing meridian over. This operation is yet another isotopy inside $S^3(L^{\ell - 1})$ 
which yields a computable presentation. Therefore, thanks to the first step, the invariant does not change. 
In the end we get the original final presentation plus some stabilizing meridian linked to the rest of the graph. All these meridians
can be slid back to their initial positions, and once again this operation is an isotopy inside $S^3(L')$.

The third step is the general case: now what we have to do is to blow-up two meridians of a component of $L$ in such a way that its 
framing does not change. Then we can consider these new curves as part of $T$, falling back into the previous case,
we can prove that we can undo our initial operation and that the result is not affected by our changes.

\subsection*{Extension to all compatible triples}

 There exist of course compatible triples which do not admit computable presentations. 
 In order to include also this case in the construction we can build a second invariant $\En^0_{\Cat}$
 which is defined for all compatible triples.
 
 \begin{remark}
  In categories where the quantum dimension of the objects of $\A$ is always zero (such as the categories in our example) 
  this second invariant will vanish on all triples which admit computable presentations.
  Therefore in this case one should continue to use $\En_{\Cat}$ to get topological informations.
 \end{remark}
 
 For the definition of $\En^0_{\Cat}$ we will need the concept of connected sum of compatible triples.
 Let $(M_1,T_1,\omega_1)$ and $(M_2,T_2,\omega_2)$ be two compatible triples, let $M_3 = M_1 \# M_2$ be the connected
 sum along balls $B_i$ inside $M_i \smallsetminus T_i$ for $i = 1,2$ and set $T_3 = T_1 \sqcup T_2$. Then we have the 
 chain of isomorphisms
 \begin{align*}
  H_1(M_3 \smallsetminus T_3) &\cong H_1(M_1 \smallsetminus (B_1 \cup T_1)) \oplus H_1(M_2 \smallsetminus (B_2 \cup T_2)) \\
  &\cong  H_1 (M_1 \smallsetminus T_1) \oplus H_1(M_2 \smallsetminus T_2)
 \end{align*}
 where the first one is induced by a Mayer-Vietoris sequence and the second
 one comes from excision. These maps induce an isomorphism 
 \[
  H^1 (M_3 \smallsetminus T_3 ; G)
  \cong  H^1 (M_1 \smallsetminus T_1 ; G) \oplus H^1 (M_2 \smallsetminus T_2 ; G).
 \]
 Finally let $\omega_3$ be the unique element of $H^1 (M_3 \smallsetminus T_3 ; G)$ which restricts to $\omega_i$ on 
 $H^1 (M_i \smallsetminus T_i ; G)$ for $i=1,2$ via the previous isomorphism. The connected sum of 
 $(M_1, T_1, \omega_1 )$ and $(M_2 , T_2 , \omega_2 )$ is defined as $(M_3 , T_3 , \omega_3 )$. 
 Now if the compatible triple $(M,T,\omega)$ does not admit any computable presentation consider the triple $(S^3,u_V,\omega_V)$
 where $u_V$ is a 0-framed unknot in $S^3$ colored by $V \in \A$ and $\omega_V$ is the unique cohomology class in
 $H^1(S^3 \smallsetminus u_V;G)$ which makes the previous triple into a compatible one. Then we can define 
 $\En^0_{\Cat}(M,T,\omega)$ to be
 \[
  \frac{\En_{\Cat}((M,T,\omega) \# (S^3,u_V,\omega_V))}{\ed(V)}.
 \]

\begin{remark}
\begin{enumerate}
 \item Just like before, when $\Cat = U^H_q(\sldue)$-mod with $q = e^{\frac{\pi i}{r}}$ we write $\En^0_r$ instead of 
 $\En^0_{U^H_q(\sldue)\mathrm{-mod}}$.
 As claimed earlier, $\En^0_r$ vanishes on computable presentations because in this category $F'$ vanishes on split $\A$-graphs, i.e.
 if $T$ and $T'$ are completely disjoint $\A$-graphs then we have $F'(T \sqcup T') = F'(T )F(T') = 0$.
 \item As it was mentioned in the abstract, $\En^0_r$ coincides with $\tau_r$ in a lot
 of cases, though in general their equality remains conjectural.
\end{enumerate}
\end{remark}

\section{Case r = 2: Alexander polynomial, Reidemeister torsion and lens spaces}
\sectionmark{Case r=2}

\begin{figure}[hbtp]
\centering
\includegraphics[width=\textwidth]{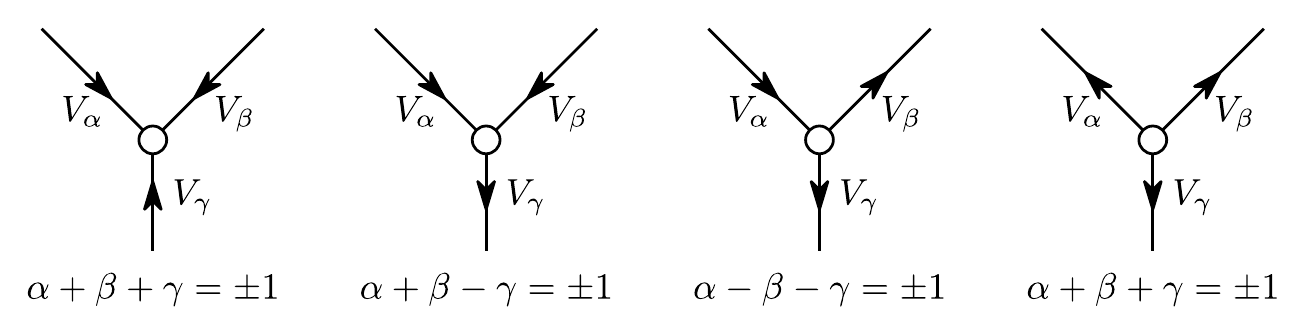}\label{trivalent}
\caption{Admissible colorings}
\end{figure}

\begin{figure}[tbp]
\centering
\includegraphics{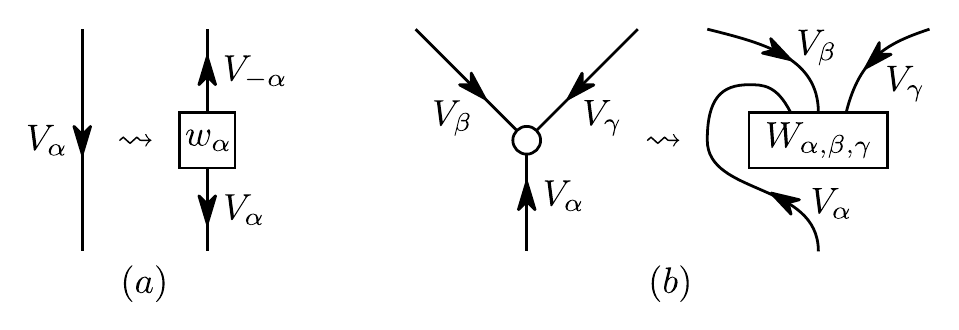}\label{sostituzioni}
\caption{$\Cat$-colored ribbon graph $T_{\Gamma}$ obtained from trivalent graph $\Gamma$}
\end{figure}

For the special case $r = 2$ we have that $q = i$ and the modified $\A$-graph invariant $F'$ associated with the category 
$U^H_i(\sldue)$-mod 
can be related to the multivariable Alexander polynomial. This fact, which was first observed by Murakami in \cite{murakami},
is exposed in detail in Viro's paper \cite{viro}. He defines an \textit{Alexander invariant} $\underline{\Delta}^2$ 
for oriented trivalent graphs equipped with the following additional structure:
\begin{enumerate}
 \item a half-integer framing (half-twists are allowed too);
 \item a coloring with typical $U^H_i(\sldue)$-modules satisfying a condition like the ones 
  shown in Figure \ref{trivalent} around each vertex;
 \item a cyclic ordering of the (germs of the) edges around each vertex.
\end{enumerate}
Viro's construction
uses a functor which is similar to the Reshetikhin-Turaev one, 
though the source category is not the category of colored ribbon graphs. 
It is indeed a category $\G^2$ whose objects are the objects of $\Rib_{U^H_i(\sldue)\mathrm{-mod}}$ 
which feature only typical colors and whose morphisms are (isotopy classes of) a non-closed version of the graphs mentioned above. 
In particular all vertices are either 3-valent (internal
vertices) or 1-valent (boundary vertices). 
If such a graph $\Gamma$ is closed, i.e. if it does not contain boundary vertices, 
and if its framing yields an orientable surface, then we can associate with it an $\A$-graph
$T_{\Gamma}$ defined as follows: 
consider the ordered basis 
$\{ v_0^{\alpha}, v_1^{\alpha} := \frac{i^{\frac{\alpha + 1}{2}}}{[\alpha + 1]} F v_0^{\alpha} \}$
of the typical module
$V_{\alpha}$ and let $\{ \varphi^0_{\alpha}, \varphi^1_{\alpha} \}$ be its dual basis in $V_{\alpha}^*$. Then 
we have an isomorphism $w_{\alpha} : V_{\alpha} \rightarrow V_{-\alpha}^*$ given by 
$v_j^{\alpha} \mapsto i^{- \left( \frac{\alpha + 1}{2} - j \right)} \varphi_{-\alpha}^{1-j}$ for $j = 0,1$.
Moreover, every time $\alpha, \beta, \gamma \in \mathbb C \smallsetminus (2 \Z + 1)$ satisfy
$\alpha + \beta + \gamma = \pm 1$, we can consider the morphism $W_{\alpha,\beta,\gamma} : \mathbb C \rightarrow V_{\alpha} \otimes
V_{\beta} \otimes V_{\gamma}$ mapping 1 to $\sum_{2(j+k-h) = \alpha + \beta + \gamma + 1} C^{\alpha,\beta,\gamma}_{j,k,h}
v^{\alpha}_j \otimes v^{\beta}_k \otimes v^{\gamma}_h$ where the coefficients $C^{\alpha,\beta,\gamma}_{j,k,h}$ are 
derived from the Clebsch-Gordan quantum coefficients (compare with \cite{costantino-murakami}) and are defined as 
\begin{align*}
 &(-1)^{k-h} i^{\frac{\beta(k-1) - \alpha(j+1) + 2(k+h-j-1) + j^2 - k^2}{2}}
 \cdot
 \left[ \begin{array}{c}
         1-\gamma \\ 1 - \gamma - h
        \end{array} \right]^{-1}
 \left[ \begin{array}{c}
         1-\gamma \\ \frac{\alpha + \beta - \gamma + 1}{2}
        \end{array} \right] \\
 &\cdot
 \sum_{t+s=h} (-1)^t i^{\frac{(2t-h)(2-\gamma-h)}{2}}
 \left[ \begin{array}{c}
         \frac{\alpha + \beta + \gamma + 1}{2} \\ j - t
        \end{array} \right] 
 \left[ \begin{array}{c}
         \alpha - j + t + 1 \\ \alpha - j + 1
        \end{array} \right]
 \left[ \begin{array}{c}
         \beta - k + s + 1 \\ \beta - k + 1
        \end{array} \right].
\end{align*}
Then we can construct $T_{\Gamma}$ by replacing each edge of $\Gamma$ which is not a connected component as shown in 
Figure \ref{sostituzioni}.(a) and each trivalent vertex of $\Gamma$ as shown in Figure \ref{sostituzioni}.(b).

\begin{proposition}\label{alexander}
 $F'(T_{\Gamma}) = (-2i)^{1 - v/2} \underline{\Delta}^2 (\overleftarrow{\Gamma})$ where $v$ is the number of vertices of $\Gamma$
 and $\overleftarrow{\Gamma}$ is obtained from $\Gamma$ by inverting the orientation on each edge.
\end{proposition}

This result is obtained by checking that the two expressions coincide for a set of elementary graphs (the trivial one, the 
$\Theta$-graph and the tetrahedron graph) and by checking that they both satisfy the same set of relations which reduce the 
computation for an arbitrary graph to elementary ones (see \cite{cgp1} and \cite{bcgp} for details).

If $L = L_1 \sqcup \ldots \sqcup L_m$ is an oriented colored framed link whose $j$-th component $L_j$ is colored with 
the typical module $V_{\alpha_j}$ then Viro shows that
\[
 \underline{\Delta}^2 (L) = \nabla_L(i^{1+\alpha_1}, \ldots, i^{1+\alpha_m}) 
 \cdot \; i^{\sum_{j,h = 1}^m \frac{\alpha_j \alpha_h - 1}{2} \ell k(L_j,L_h)},
\]
where $\nabla_L$ is the Alexander-Conway function of $L$.
Therefore, if the framing is integral, Proposition \ref{alexander} immediately gives
\[
 F'(L) = (-2i) \nabla_L(i^{1-\alpha_1}, \ldots, i^{1-\alpha_m}) 
 \cdot \; i^{\sum_{j,h = 1}^m \frac{\alpha_j \alpha_h - 1}{2} \ell k(L_j,L_h)}.
\]

Now since $\Cat_{\bar{1}}$ is semisimple inside $U^H_i(\sldue)$-mod we can take the critical set $X \subset \mathbb C/2\Z$
to be just $\{ \bar{0} \}$. Therefore, thanks to Remark \ref{computable}, every triple of the form $(M,\varnothing,\omega)$ with
$\omega \neq 0$ is compatible and admits a computable presentation. In particular some computation (compare with \cite{bcgp}) 
yields
\begin{align*}
 \En_2(M,\varnothing,\omega) = 2 &\cdot 4^{m-\sigma_+(L)-\sigma_-(L)} 
 \cdot i^{\sigma_-(L)-\sigma_+(L)-m-1} \cdot \left( \prod_{j=1}^m \frac{1}{i^{\alpha_j}-i^{-\alpha_j}}\right) \\
 &\cdot \nabla_L(i^{\alpha_1},\ldots,i^{\alpha_m}) \cdot i^{\sum_{j,h = 1}^m \frac{\alpha_j (\alpha_h + 2)}{2} \ell k(L_j,L_h)},
\end{align*}
where $L = L_1 \sqcup \ldots \sqcup L_m$ is a surgery presentation for $M$ and $\alpha_j := \langle \omega, \mu_j \rangle$.
Thus $\En_2$ recovers the Alexander-Conway function, which is known to be related to the Reidemeister torsion. Moreover
$\En_2$ yields a canonical normalization of the Reidemeister torsion which
fixes the scalar indeterminacy. Indeed recall that the refined abelian Reidemeister torsion of $M$ defined by Turaev (see 
\cite{turaev-torsion} for example) is determined by the choice of a homomorphism $\varphi : H_1(M) \rightarrow \mathbb C^*$,
of a homology orientation $\omega_M$ for $M$ and of a $\Spin^c$-structure $\sigma \in \Spin^c(M)$ 
(or equivalently of an Euler structure on $M$).
We write $\tau^{\varphi}(M,\omega_M,\sigma)$ or, if $M$ is oriented and we pick the canonical homology orientation associated with
the orientation of $M$, simply $\tau^{\varphi}(M,\sigma)$.
Now, if $(M,\varnothing,\omega)$ is a compatible triple as above, we can use the non-zero cohomology class $\omega$ to define the 
homomorphism $\varphi_{\omega} : H_1(M) \rightarrow \mathbb C^*$ given by $h \mapsto e^{i \pi \langle \omega, h \rangle} = 
i^{2 \langle \omega, h \rangle}$.

\begin{theorem}
 Let $M$ be a closed oriented 3-manifold endowed with a non-trivial
 cohomology class $\omega \in H^1(M;\mathbb C/2\Z)$. Then for any complex spin structure $\sigma$ in $\mathrm{Spin}^c(M)$ we have
 \[
  \tau^{\varphi_{\omega}}(M,\sigma) = \frac{i^{b_1(M) + 4 \psi_{M,\sigma}(\omega) + 1 }}{2 \cdot 4^{b_1(M)}}
  \cdot \En_2(M,\varnothing,\omega),
 \]
 where 
 \[
  \psi_{M,\sigma} : H^1(M;\mathbb C/2\Z) \rightarrow \mathbb C/\Z
 \]
 is the homomorphism obtained by first extending
 Deloup and Massuyeau's quadratic linking function 
 \[
  \varphi_{M,\sigma} : H_2(M;\bbQ/\Z) \rightarrow \bbQ/\Z
 \]
 (compare with 
 \cite{deloup-massuyeau}, Definition 2.2)
 to a homomorphism 
 \[
  \varphi^{\mathbb C}_{M,\sigma} : H_2(M;\mathbb C/\Z) \rightarrow \mathbb C/\Z,
 \]
 and then by considering
 the composition $\varphi^{\mathbb C}_{M,\bar{\sigma}} \circ \frac 12 \circ D$ where 
 \[
  D : H_2(M;\mathbb C/2\Z)
 \rightarrow H^1(M;\mathbb C/2\Z)
 \]
 is Poincar\'{e} duality, 
 \[
  \frac 12 : H^1(M;\mathbb C/2\Z) \rightarrow H^1(M;\mathbb C/\Z)
 \]
 is induced by the \textquotedblleft division by 2\textquotedblright \hspace{1pt}
 isomorphism between $\mathbb C/2 \Z$ and $\mathbb C/\Z$ and $\bar{\sigma}$ is the image of $\sigma$
 under the standard involution of $\Spin^c(M)$
\end{theorem}

This is proven by using the surgery formula for the Reidemeister torsion (see \cite{turaev-torsion}, section VIII.2, equation (2.b)):
if $L$ is a computable surgery presentation for $(M,\varnothing,\omega)$
then
\[
 \tau^{\varphi_{\omega}}(M,\sigma) = (-1)^{2m - \sigma_+(L)} 
 \cdot \prod_{j=1}^m 
 \frac{i^{\alpha_j(k_j - 1)}}{i^{\alpha_j} - i^{-\alpha_j}} \cdot \nabla_L(i^{\alpha_1},\ldots,i^{\alpha_m}),
\]
where $\alpha_j := \langle \omega,\mu_j \rangle$ and $k_1, \ldots, k_m$ are the charges of the $\Spin^c$-structure $\sigma$
(see \cite{turaev-torsion}, section VII.2.2 for a definition).

We conclude with a Proposition which gives the value of $\En_2$ for lens spaces and can be used to follow the path of 
the classical proof of their classification.

\begin{proposition}
 Let $p > q > 0$ be two coprime integers, let $L(p,q)$ be a lens space and consider a non zero cohomology class 
 $\omega \in H^1(L(p,q);\mathbb C/2\Z)$. Then
 \[
  \En_2(L(p,q),\omega) = \frac{(-1)^{k(\omega)} e^{\frac{i \pi k(\omega)^2 p}{q}}}{2i \sin\frac{\pi k(\omega) q}{p}
  \sin\frac{\pi k(\omega)}{p}}
 \]
 for some $k(\omega) \in \Z \smallsetminus p\Z$.
\end{proposition}

\begin{corollary}
 $\En_2$ classifies lens spaces.
\end{corollary}

\end{document}